\documentclass[a4paper,12pt,reqno]{amsart}
\usepackage{amssymb}
\usepackage{amsmath}

\def\i{\lrcorner}

\def\.{\cdot}

\def\n{\nabla}

\def\beq{\begin{equation}}
\def\eeq{\end{equation}}
\def\bea{\begin{eqnarray*}}
\def\eea{\end{eqnarray*}}
\def\x{\times}

\def\L{\Lambda}

\def\r{\end{proof}}

\def\SO{\rm SO}
\def\SU{\rm SU}
\def\Sp{\rm Sp}
\def\U{\rm U}

\def\d{{\delta}}
\def\dm{d_1}
\def\dem{\d_1}
\def\dn{d_2}
\def\den{\d_2}

\def\es{\,\lrcorner\;}

\def\be{\begin{equation}}

\def\ee{\end{equation}}
\def\Lie{{\mathcal L\,}}

\def\Hol{\mathrm{Hol}}

\def\K{{\mathfrak K}}
\def\P{{\mathfrak P}}
\def\T{{\mathfrak T}}

\newtheorem{ede}{Definition}[section]

\newtheorem{epr}[ede]{Proposition}

\newtheorem{ath}[ede]{Theorem}

\newtheorem{ecor}[ede]{Corollary}
 
\title{Twistor Forms on Riemannian Products}
\author{Andrei Moroianu and Uwe Semmelmann}
\address{Andrei Moroianu \\ CMLS\\ {\'E}cole Polytechnique \\ UMR 7640 du CNRS
\\ 91128 Palaiseau \\ France}
\email{am@math.polytechnique.fr}

\address{Uwe Semmelmann\\ Universit\"at zu K\"oln,
Mathematisches Institut \\ Weyertal 86--90
D--50931 K\"oln, Germany}
\email{uwe.semmelmann@math.uni-koeln.de}

\begin{document}

\begin{abstract}
We study twistor forms on products of compact Riemannian
manifolds and show that they are defined by Killing forms on 
the factors. The main result of this note is a necessary step in the
classification of compact Riemannian manifolds with non-generic holonomy
carrying twistor forms.

\vspace{0.2cm}


\end{abstract}

\maketitle


\section{Twistor Forms on Riemannian Manifolds}

A {\em twistor $p$-form} on a Riemannian manifold $\,(M^n,\,g)$ is a
smooth section $\psi$ of $\Lambda^pT^*M$ whose covariant derivative
only depends on its differential $d\psi$ and codifferential $\d\psi$. 
More precisely, $\psi$ satisfies the equation
\begin{equation}\label{killing}
\nabla_X\,\psi\;=\;
\tfrac{1}{  p+1}X\es d\psi -
\tfrac{1}{  n-p+1} X^\flat\wedge\d \psi,
\end{equation}
for all vector fields $X$, where $X^\flat$ denotes the metric dual of
$X$.

If the $p$-form $\psi$ is in addition coclosed ({\em i.e.}
$\d\psi=0$), then it is called a {\em Killing $p$-form}. We denote by
$\T(M)$, $\K(M)$ and $\P(M)$ the spaces of twistor, Killing and
parallel forms on $M$ respectively. Notice that $\T(M)$ is preserved
by Hodge duality, and that the 
Hodge dual of a Killing form is a closed twistor form.
For a comprehensive introduction to twistor forms,
see \cite{uwe}.  

A few years ago, a program of classification of twistor forms on
compact manifolds was started. By the de Rham decomposition theorem, every
simply connected Riemannian manifold is a Riemannian product of
{\em irreducible} manifolds. Moreover, the Berger-Simons holonomy theorem
(see \cite{besse}, p. 300)
implies that any simply connected irreducible Riemannian manifold is
either symmetric or has holonomy $\SO_n$, $\U_m$, $\SU_m$, $\Sp_k$,
$Sp_k\.\Sp_1$, $G_2$ or $Spin_7$. Killing forms on symmetric spaces
were studied in \cite{bms}. Twistor forms on K\"ahler manifolds
(covering the holonomies $\U_m$, $\SU_m$, and $\Sp_k$) were
described in \cite{ms1}, and Killing forms on quaternion-K\"ahler
manifolds (holonomy $Sp_k\.\Sp_1$) or Joyce manifolds (holonomies $G_2$
or $Spin_7$) were studied in \cite{ms2} and \cite{uwe1}
respectively. In Theorem \ref{th} below, we prove that the general
case (twistor forms on a Riemannian product of compact manifolds)
reduces to the study of Killing 
forms on the factors. By the discussion above, besides the case of
generic holonomy ($\SO_n$), all other cases are fully understood.


\section{The Main Result}

Let $M=M_1\x M_2$ be the Riemannian product of two compact Riemannian manifolds
$(M_1,g_1)$ and $(M_2,g_2)$ of dimensions $m$ and $n$
respectively. We denote by $\pi_i$ the projection $\pi_i:M\to
M_i$. From (\ref{killing}) it is clear that $\pi_i ^*(\K(M_i))\subset
\K(M)$, so the space 
$$\K_0(M):=\pi_1 ^*(\K(M_1))+\pi_2 ^*(\K(M_2))+\P(M)$$
is a subspace of $\K(M)$. For later use, we give the following
description of $\pi_i ^*(\K(M_i))$:
\begin{equation}\label{c1}\pi_1 ^*(\K(M_1))=\{u\in\K(M)\ |\ \n_X u=0,\
  \forall X\in TM_2\}\end{equation}
and 
 \begin{equation}\label{c2}\pi_2 ^*(\K(M_2))=\{u\in\K(M)\ |\ \n_X
   u=0,\ \forall X\in TM_1\}.\end{equation}

The aim of this 
note is to prove the following result:

\begin{ath} \label{th} Every twistor form on $M$
  is a sum of forms of the following types: parallel forms,
  pull-backs of Killing forms on $M_1$ or $M_2$, and Hodge duals of
  them. In other words, $\T(M)=\K_0(M)+*\K_0(M).$
\end{ath}

\begin{proof} Since $\K_0(M)\subset \K(M)\subset \T(M)$ and
  $*\T(M)=\T(M)$, we clearly have $\K_0(M)+*\K_0(M)\subset \T(M)$. It
  remains to prove the reverse inclusion.
Let us define the differential operators
$$\dm=\sum_{i=1}^m e_i ^\flat\wedge\n_{e_i},\qquad\qquad 
\dn=\sum_{j=1}^n f_j ^\flat\wedge\n_{f_j},$$
where $\{e_i\}$ and $\{f_j\}$ denote local orthonormal basis of the
tangent distributions to $M_1$ and $M_2$. Using the Fubini
theorem, we easily see that the adjoint operators to $\dm$ and $\dn$ are
$$\dem=-\sum_{i=1}^m e_i\i\n_{e_i},\qquad\qquad 
\den=-\sum_{j=1}^n f_j\i\n_{f_j}.$$
The following relations are straightforward:
$$d^M=\dm+\dn,\qquad \d^M=\dem+\den,\qquad
(\dm)^2=(\dn)^2=(\dem)^2=(\den)^2=0,$$
$$0=\dm\dn+\dn\dm=\dem\den+\den\dem,\qquad 0=\dm\den+\den\dm=\dem\dn+\dn\dem.$$
The vector bundle $\L^pM$ decomposes naturally as
$$\L^pM\cong\oplus_{i=0}^p \L^{i,p-i}M,$$
where $\L^{i,p-i}M\cong\L^iM_1\otimes\L^{p-i}M_2.$
Obviously, $\dm$ and $\dem$ map $\L^{i,p-i}M$ to $\L^{i+1,p-i}M$ and
$\L^{i-1,p-i}M$ respectively, and $\dn$ and $\den$ map $\L^{i,p-i}M$
to $\L^{i,p-i+1}M$ and $\L^{i,p-i-1}M$ respectively.

With respect to the above decomposition, every $p$-form can be written
$u=u_0+\ldots +u_p$, where $u_i\in \L^iM_1\otimes\L^{p-i}M_2$. From
now on, $u$ will denote a twistor $p$-form $u\in\T(M)$,
with $1\le p\le n+m-1$. The twistor equation reads
\begin{equation}\label{tw}
\n_Xu=\frac{1}{p+1}X\i(\dm u+\dn u)-\frac{1}{m+n-p+1}X\wedge(\dem
u+\den u),\quad \forall X\in TM.
\end{equation}
By projection onto the different irreducible
components of $\L^pM$, (\ref{tw}) can be translated into
the following two systems of equations:
\begin{equation}\label{tw1}
\n_Xu_k=\frac{1}{p+1}X\i(\dm u_k+\dn u_{k+1})-\frac{1}{m+n-p+1}X\wedge(\dem
u_k+\den u_{k-1}),\quad \forall X \in TM_1,
\end{equation}
and
\begin{equation}\label{tw2}
\n_Xu_k=\frac{1}{p+1}X\i(\dm u_{k-1}+\dn u_{k})-\frac{1}{m+n-p+1}X\wedge(\dem
u_{k+1}+\den u_{k}),\quad \forall X \in TM_2.
\end{equation}
Recall that if $u$ is any $k$-form and $\{e_1,\ldots,e_m\}$ is an
orthonormal basis on a manifold $M$, then 
\begin{equation}\label{c}
\sum_{i=1}^m e_i ^\flat\wedge e_i\i\omega=k\omega.
\end{equation}
Taking the wedge product with $X^\flat$ in (\ref{tw1}) and summing over an
orthonormal basis of $TM_1$ yields 
$$\dm u_k=\sum_{i=1}^m
  e_i\wedge\n_{e_i}u_k=\frac{1}{p+1}\sum_{i=1}^m e_i\wedge e_i\i(\dm
  u_k+\dn u_{k+1}) 
\stackrel{(\ref{c})}{=}\frac{k+1}{p+1}(\dm u_k+\dn
u_{k+1})
$$
 so 
\begin{equation}\label{t1}
(p-k)\dm u_k=(k+1)\dn u_{k+1}.
\end{equation}
Similarly, taking the interior product with $X$ and summing over an
orthonormal basis of $TM_1$ yields $\dem u_k=\frac{m-k+1}{m+n-p+1}(\dem
u_k+\den u_{k-1})$, thus
\begin{equation}\label{t2}
(n+k-p)\dem u_k=(m-k+1)\den u_{k-1}.
\end{equation}
We distinguish three cases:

{\bf Case I.} Suppose that $p$ is strictly smaller than $m$ and $n$. For $k<p$,
(\ref{t1}) and (\ref{t2}) imply
\begin{equation}\label{tt}
\dem\dm u_k=\frac{k+1}{p-k}\dem\dn u_{k+1}=-\frac{k+1}{p-k}\dn\dem
u_{k+1}
=-\frac{(k+1)(m-k)}{(p-k)(n+k-p+1)}\dn\den u_k.\end{equation}
Integrating over $M$ yields $0=\dm u_k=\den u_k,\ \forall k<p$. Similarly one
gets $0=\dn u_k=\dem u_k,\ \forall k>0$. Moreover, we have $0=\den u_p
=\dem u_0$ (tautologically), so in particular $\dem u_k=\den u_k=0,\
\forall k$. From (\ref{tw1}) and (\ref{tw2}), together with (\ref{c1})
and (\ref{c2}), we see that $u_1,\ldots,
u_{p-1}\in\P(M)$, $u_0\in\pi_2^*(\K(M_2))$ and
$u_p\in\pi_1^*(\K(M_1))$, so $u\in\K_0(M)$. 

{\bf Case II.} Suppose that $p$ is strictly larger than $m$ and
$n$. Since the Hodge dual $*u$ of $u$ is a twistor
$(m+n-p)$-form and $m+n-p$ is strictly smaller than $m$ and $n$,
the first case implies that $*u\in\K_0(M)$, so $u\in *\K_0(M)$.

{\bf Case III.} If $p$ is a number between $m$ and 
$n$, we may suppose without loss of generality that $m\le p\le
n$. Obviously $u_{m+1}=\ldots =u_p$=0. Using
(\ref{tt}) and integrating over $M$, we
obtain  that $0=\dm u_k=\den u_k$ for $0\le k\le m-1$ and similarly,
$0=\dn u_k=\dem 
u_k$ for $1\le k\le m$. As before, (\ref{tw1}) and (\ref{tw2}),
together with (\ref{c1}) and (\ref{c2}),  show
that $u_1,\ldots,
u_{m-1}\in\P(M)$, $u_0\in\pi_2^*(\K(M_2))$, and
$*u_m\in\pi_2^*(\K(M_2))$. This proves the theorem. 

\r

As an application of this result, we have the following:

\begin{epr}\label{pr} Let $(M^n,g)$ be a compact simply connected Riemannian
  manifold. If $M$ 
  carries a conformal vector field which is not Killing, then
  $\Hol(M)=\SO_n$. 
\end{epr}
\begin{proof} Assume first that $(M,g)=(M_1,g_1)\times (M_2,g_2)$ is a
Riemannian product with $\dim(M_1),\ \dim(M_2)\ge1.$ Then, taking into
account that the isomorphism between $1$-forms and vector fields defined by the
Riemannian metric maps twistor forms to conformal vector fields and
Killing forms to Killing vector fields, Theorem \ref{th} implies that
every conformal vector field on $M$ is a Killing vector field. Thus
$M$ is irreducible. 

Assume next that $\Hol(M)\ne\SO_n$. From the Berger-Simons holonomy
theorem (\cite{besse}, p. 300), $M$ is either an irreducible symmetric
space (in particular 
Einstein), or its holonomy group is $\U_m$, $\SU_m$, $\Sp_k$,
$Sp_k\.\Sp_1$, $G_2$ or $Spin_7$. In the first three cases the manifold is
K\"ahler and in the last three cases it is Einstein.
Now, two classical results state that a conformal vector field on a
compact manifold $M$ is 
already a Killing vector field if $M$ is K\"ahler (see \cite{li},
p. 148) or if $M$ is 
Einstein and not isometric to the round sphere (see \cite{na}, \cite{na-ya}).

The only possibility left is therefore $\Hol(M)=\SO_n$.

\r

{\bf Example.} Take any compact simply connected Riemannian
  manifold $(M^n,g)$ carrying a Killing vector field $\xi$ and let $f$
  be a function on $M$ such that $\xi(f)$ is not identically
  zero. Since $\Lie_\xi(e ^{2f} g)=2\xi(f)e
  ^{2f} g$, $\xi$ is a conformal vector field on
  $(M,e ^{2f} g)$ which is not Killing. From Proposition \ref{pr},
  $(M,e ^{2f} g)$ has holonomy $\SO_n$.    

\begin{ecor}\label{cor1}
Let $(M^n,g)$ be a compact simply connected homogeneous Riemannian
  manifold. Then for every non-constant function $f$ on $M$,
  $(M,e ^{2f} g)$ has holonomy $\SO_n$.  
\end{ecor}
\begin{proof} Since $f$ is non-constant, there exists $x\in M$ such
  that $df_x\ne 0$. Killing vector fields on $M$ span the tangent
  spaces at each point, so in particular there exist a Killing vector field
$\xi$ such that $\xi(f)$ is not identically
  zero. The corollary then follows from the example above.
\r

 \labelsep .5cm

\end{document}